\documentclass[twoside]{amsart}
\usepackage{amssymb,amsmath,amscd}
\usepackage{amsthm}
\usepackage{amsfonts}
\usepackage{enumerate}
\usepackage{amsbsy}
\usepackage{geometry}
\usepackage[utf8]{inputenc}
\usepackage[T1]{fontenc}
\usepackage{microtype}

\newtheorem{theorem}{Theorem}[section]

\newtheorem{definition}[theorem]{Definition}

\newtheorem{problem}[theorem]{Problem}
\numberwithin{equation}{section}
\usepackage{caption} % For caption spacing
\usepackage{subcaption} % For sub-figures
\usepackage{graphicx}
\usepackage{pgfplots}
\usepackage[all]{nowidow}
\usepackage[utf8]{inputenc}
\usepackage{tikz}
\usepackage[inline]{enumitem}
\usepackage{blindtext}
\usepackage{amssymb}		
\usepackage{ulem}
\usepackage{amsfonts}

\usepackage[hidelinks,colorlinks=true,linkcolor=blue,citecolor=red, urlcolor=blue]{hyperref}
\numberwithin{equation}{section}
\pgfplotsset{compat=1.18}

\begin{document}
\title{An Analytic Construction of Random Variables in Lebesgue Spaces}

\author[Hugo Reyna-Castañeda]{Hugo Guadalupe Reyna-Castañeda}
\address{Departamento de Matemáticas, Facultad de Ciencias, Universidad Nacional Autónoma de México, Mexico City, Mexico}
\curraddr{}
\email{hugoreyna46@ciencias.unam.mx}
\thanks{}

%    author two information
\author{Mar\'ia de los \'Angeles Sandoval-Romero}
\address{Departamento de Matemáticas, Facultad de Ciencias, Universidad Nacional Autónoma de México, Mexico City, Mexico}
\curraddr{}
\email{selegna@ciencias.unam.mx}
\thanks{}

\begin{abstract}
    This work develops, from an analytic perspective, the construction of the concept of a random variable in Lebesgue spaces $L^{p}$. Starting from an elementary problem, it motivates the need to extend the classical notions of probability —measurability, integrability, and expectation— to the context of functions taking values in $L^{p}$. To this end, the Lebesgue spaces and their fundamental properties are introduced, emphasizing Pettis’s measurability theorem as a key tool for characterizing $L^{p}$--valued functions. Subsequently, the expectation of such random variables is defined through a functional formulation of the Riesz representation theorem, allowing the Bochner integral to be interpreted as a natural generalization of the classical expectation. Altogether, the text offers a clear and coherent framework for understanding how the spaces $L^{p}$ provide an analytic setting in which probability naturally intertwines with the linear structure of functional analysis.
\end{abstract}

\maketitle
\setlength{\parskip}{0.25em} 

\section{Introduction}

Let us begin with a simple problem in the context of the probability space $\big((0,1),\mathcal{L}(0,1),\lambda\big)$, where $\mathcal{L}(0,1)$ denotes the Lebesgue $\sigma$--algebra and $\lambda$ the Lebesgue measure.

\begin{problem} \label{Problem1}
Let $\varepsilon\in(0,1)$ and $p\in(1,\infty)$ be fixed.  
For which values of $\omega\in(0,1)$ is it possible to compute the probability that the area under the curve given by the graph of the function $\mathbf{1}_{(0,\omega)}$ is less than $\varepsilon^{p}$?
\end{problem}

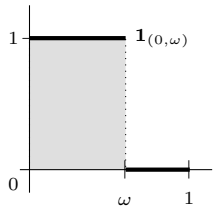
\begin{figure}[!ht]
    \centering
	\begin{tikzpicture}[xscale=0.85,yscale=0.875]
    %\draw(4.5,1.25)node[above]{$_{\chi}$};
    \draw [fill=gray!25,  dotted]  (6.5,0)--(6.5,2)--(8,2)--(8,0);
	\draw[-] (6.35,0)--(9.35,0);
	\draw[-] (6.5,-0.35)--(6.5,2.5);
\draw (6.5,2) node{$_{-}$}; \draw (6.5,2) node[left] {$_{1}$}; \draw (6.5,-0.22) node[left] {$_{0}$};
\draw (8,0) node{$_{|}$}; \draw (8,-0.27) node[below] {$_{\omega}$};
\draw (9,0) node{$_{|}$}; \draw (9,-0.2) node[below] {$_{1}$};
\draw [ultra thick,black] (8,0)--(9,0);
\draw [ultra thick,black] (6.5,2)--(8,2);
\draw (8,2) node[right]{$_{\mathbf{1}_{(0,\omega)}}$};
\end{tikzpicture}
%\caption{Representación gráfica}
\end{figure}

We can reformulate this question by observing that the area under the curve corresponds to the integral of the function under consideration. Thus, the problem can be expressed as follows:

\begin{problem}\label{Problem2}
Let $\varepsilon\in(0,1)$ and $p\in(1,\infty)$ be fixed.  
Can we compute the probability of the set
$$
\left\{\,\omega\in(0,1)\;:\;\int_{0}^{1}\mathbf{1}_{(0,\omega)}(t)\,dt<\varepsilon^{p}\right\}?
$$
\end{problem}

In an elementary probability course, problems of this type are usually treated through the notion of a random variable. Recall that, given a probability space $(\Omega,\mathcal{F},\mathbb{P})$, a function $X:\Omega\to\mathbb{R}$ is called a random variable if $X^{-1}(B)\in\mathcal{F}$ for every Borel set $B\in\mathcal{B}(\mathbb{R})$.  

This measurability requirement ensures that the probability of events defined through $X$ is well determined, and it naturally leads to the result that allows the measure $\mathbb{P}$ to be transported from $\Omega$ to $\mathbb{R}$ (see \cite[Theorem~8.5]{JacodProtter2004}).

\begin{theorem} \label{Medida-Real-Distribucion}
Let $(\Omega,\mathcal{F},\mathbb{P})$ be a probability space and $X:\Omega\to\mathbb{R}$ a random variable.  
Then the mapping 
$$
X_{\sharp}\mathbb{P}(B):=\mathbb{P}(X^{-1}(B)), \qquad B\in\mathcal{B}(\mathbb{R}),
$$
defines a probability measure on $(\mathbb{R},\mathcal{B}(\mathbb{R}))$.
\end{theorem}

Let us now return to our initial problem.  
Consider the set
$$
\mathcal{I}(0,1):=\{\,\mathbf{1}_{A}:A\in\mathcal{L}(0,1)\,\},
$$
consisting of the indicator functions of measurable subsets of $(0,1)$.  
We define the mapping $\chi:(0,1)\to\mathcal{I}(0,1)$ by
$$
\chi(\omega):=\mathbf{1}_{(0,\omega)}.
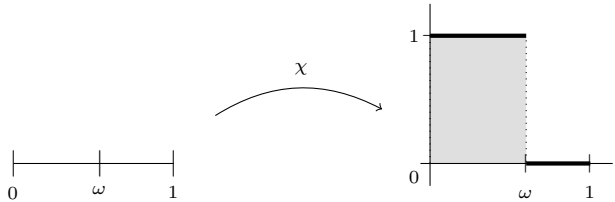
$$

\begin{figure}[!ht]
    \centering
	\begin{tikzpicture}[xscale=0.85,yscale=0.85]
	\draw[-] (0,0)--(2.5,0);
	\draw (0,0) node{$|$}; \draw (0,-0.22) node[below]{$_{0}$};
	\draw (2.5,0) node{$|$}; \draw (2.5,-0.2) node[below]{$_{1}$};
	\draw (1.35,0) node{$|$}; \draw (1.35,-0.2) node[below]{$_{\omega}$};
	%\draw[->] (0,-2)--(0,2);
	%\draw[-, densely dashed] (0,1)--(2.35,1);
%\draw[-, densely dashed] (0,-1)--(2.35,-1);
%\draw[domain= 0:2.35, thick, blue] plot(\x,{\x^(1/2)} );

    \draw[<-] (5.75,0.85) to[bend right] (3.15,0.75);
    \draw(4.5,1.25)node[above]{$_{\chi}$};
    \draw [fill=gray!25,  dotted]  (6.5,0)--(6.5,2)--(8,2)--(8,0);
	\draw[-] (6.35,0)--(9.35,0);
	\draw[-] (6.5,-0.35)--(6.5,2.5);
\draw (6.5,2) node{$_{-}$}; \draw (6.5,2) node[left] {$_{1}$}; \draw (6.5,-0.22) node[left] {$_{0}$};
\draw (8,0) node{$_{|}$}; \draw (8,-0.27) node[below] {$_{\omega}$};
\draw (9,0) node{$_{|}$}; \draw (9,-0.2) node[below] {$_{1}$};
\draw [ultra thick,black] (8,0)--(9,0);
\draw [ultra thick,black] (6.5,2)--(8,2);

\end{tikzpicture}
%\caption{Representación gráfica}
\end{figure}

Naturally, the question arises whether it is possible to apply the previous theorem to $\chi$ and, consequently, to assign a probabilistic meaning to the expression
$$
\lambda\!\left(\left\{\omega\in(0,1):\left(\int_{0}^{1}\chi(\omega)(t)\,dt\right)^{\frac{1}{p}}<\varepsilon\right\}\right).
$$

However, the answer is not immediate: the function $\chi$ does not take real values but rather functions, and therefore lies outside the classical framework in which random variables are measurable real-valued mappings. This fact, rather than representing an obstacle, reveals a conceptual opportunity. It invites us to extend the notion of a random variable to settings where the realizations are not numbers but functional objects —elements of an infinite-dimensional Banach space.

The aim of this work is precisely to develop such an extension. Building on the preceding example, we construct the concept of a random variable taking values in Lebesgue spaces, interpret measurability in the Pettis–Bochner sense, and show how the theory of vector measures allows one to define both its distribution and its expectation.

In this way, the initial problem serves as a guiding thread for introducing probability in spaces of functions, where randomness acquires a genuinely geometric dimension.

\section{From Indicator Functions to Lebesgue Spaces}

A large part of the classical concept of a real random variable is closely tied to the vector space structure of $\mathbb{R}$; for example, the sum of two real random variables is again a real random variable.

The function $\chi$ defined above takes values in the set $\mathcal{I}(0,1)$ consisting of the indicator functions of Lebesgue–measurable subsets of $(0,1)$. However, this set is not suitable if we wish to extend the usual algebraic operations from the real case, since the sum of two indicator functions is, in general, not an indicator function. Consequently, $\mathcal{I}(0,1)$ lacks a linear structure.

Nevertheless, not all is lost: the functions in $\mathcal{I}(0,1)$ share a key property, they are integrable in the sense of Lebesgue and, in fact, by the construction of Problem \ref{Problem2}, all of them belong to the space $L^{p}(0,1)$ for $p\in(1,\infty)$. This space can be interpreted as the collection of random variables $X:(0,1)\to\mathbb{R}$ whose $p$-th power is integrable, or equivalently, those for which the expectation $\mathbb{E}(|X|^{p})$ is finite.

Thus, we may naturally consider an extension of the function $\chi:(0,1)\to L^{p}(0,1)$,
$$
\chi(\omega)=\mathbf{1}_{(0,\omega)},
$$
and ask whether it is possible to generalize the notion of a random variable to functions taking values in a Lebesgue space.

To address this question, we shall work in a general setting.  
Let $(\Omega,\mathcal{F},\mathbb{P})$ be a complete probability space.

From now on, we identify two random variables $X,Y:\Omega\to\mathbb{R}$ whenever they are equal almost surely, that is,
$$
\mathbb{P}\big(\{\omega\in\Omega : X(\omega)=Y(\omega)\}\big)=1.
$$

For $p\in[1,\infty)$ we define
$$
L^{p}(\Omega):=\left\{X:\Omega\to\mathbb{R}\,:\,\int_{\Omega}|X(\upsilon)|^{p}\,d\mathbb{P}(\upsilon)<\infty\right\}.
$$

This set is a real vector space under the usual operations of addition and scalar multiplication, and it is called the \textbf{Lebesgue space} $L^{p}(\Omega)$.

To each element $X\in L^{p}(\Omega)$ we associate its $p$-th norm,
$$
\Vert X\Vert_{p}:=\left(\int_{\Omega}|X(\upsilon)|^{p}\,d\mathbb{P}(\upsilon)\right)^{1/p},
$$
with respect to which $L^{p}(\Omega)$ acquires a fundamental analytic structure (see \cite[Theorem 4.8]{Brezis}).

\begin{theorem}[Riesz--Fischer]\label{Theo-Riesz-Fischer}
The Lebesgue space $L^{p}(\Omega)$, endowed with the norm $\Vert\cdot\Vert_{p}$, is a Banach space; that is, every Cauchy sequence in $L^{p}(\Omega)$ converges to an element of the same space.
\end{theorem}

Among the most important properties of Lebesgue spaces stands out Hölder’s inequality (see \cite[Theorem 4.6]{Brezis}), which expresses the compatibility between the norms $p$ and $q$ associated with conjugate exponents.

\begin{theorem}[Hölder’s Inequality]\label{HolderInequality}
Let $p,q\in(1,\infty)$ satisfy $\frac{1}{p}+\frac{1}{q}=1$.  
If $X\in L^{p}(\Omega)$ and $Y\in L^{q}(\Omega)$, then $XY\in L^{1}(\Omega)$ and
$$
\|XY\|_{1}\leq\|X\|_{p}\,\|Y\|_{q}.
$$
\end{theorem}

This inequality, in addition to being a cornerstone of analysis, reveals a deep relationship between the spaces $L^{p}(\Omega)$ and $L^{q}(\Omega)$: every random variable $Y\in L^{q}(\Omega)$ naturally defines a linear and continuous functional on $L^{p}(\Omega)$ through the correspondence
$$
L_Y(X):=\int_{\Omega}X(\upsilon)Y(\upsilon)\,d\mathbb{P}(\upsilon).
$$

A natural question then arises: can every linear and continuous functional on $L^{p}(\Omega)$ be represented in this way?

The answer is affirmative and is summarized in the following classical result (see \cite[Theorem 4.11]{Brezis}).

\begin{theorem}[Riesz Representation Theorem]\label{The-Representation-Riesz}
Let $p,q\in(1,\infty)$ satisfy $\frac{1}{p}+\frac{1}{q}=1$, and let $L:L^{p}(\Omega)\to\mathbb{R}$ be a linear and continuous functional.  
Then there exists a unique function $Y\in L^{q}(\Omega)$ such that
$$
L(X)=\int_{\Omega}X(\upsilon)Y(\upsilon)\,d\mathbb{P}(\upsilon), \qquad \forall\,X\in L^{p}(\Omega).
$$
\end{theorem}

This theorem establishes an elegant correspondence between the continuous linear functionals on $L^{p}(\Omega)$ and the elements of $L^{q}(\Omega)$, showing that the duality between these spaces is expressed through integration.  
Thus, to each linear and continuous functional $L:L^{p}(\Omega)\to\mathbb{R}$ one can associate the unique function $Y_{L}\in L^{q}(\Omega)$ that represents it.

This idea will be essential later on: both the notion of expectation and the extension of the concept of a random variable to vector-valued contexts can be understood as manifestations of this functional relationship.

An additional structural feature of Lebesgue spaces, depending on certain conditions on the underlying probability space, is \textbf{separability}, that is, the existence of a countable dense subset (see \cite[Theorem 4.13]{Brezis}).

\begin{definition}\label{Def-Sigma-Gen}
A probability space $(\Omega,\mathcal{F},\mathbb{P})$ is said to be \textbf{$\sigma$--generated} if there exists a countable collection $\mathcal{A}=\{A_{j}\in\mathcal{F}:j\in\mathbb{N}\}$ such that $\mathcal{F}=\sigma(\mathcal{A})$.
\end{definition}

For instance, on any open interval $(a,b)\subset\mathbb{R}$, the Borel $\sigma$--algebra $\mathcal{B}(a,b)$ is $\sigma$--generated by the collection
$$
\left\{\left(a+\frac{1}{j},\,b-\frac{1}{j}\right):j\in\mathbb{N}\right\}.
$$

\begin{theorem}\label{Theo-Lp-Sep}
If $(\Omega,\mathcal{F},\mathbb{P})$ is a $\sigma$--generated probability space, then $L^{p}(\Omega)$ is separable for every $p\in(1,\infty)$.
\end{theorem}

Hence, the probability space $((0,1),\mathcal{L}(0,1),\lambda)$ is $\sigma$--generated, and therefore $L^{p}(0,1)$ is separable for all $p\in(1,\infty)$.

This property will be crucial, since separability makes it possible to construct a countable sequence of linear and continuous functionals that characterize the norm of any element in $L^{p}(\Omega)$ —a feature that may be viewed as a countable manifestation of the Hahn-Banach theorem (see \cite[Corollary 1.9.8 and Theorem 1.10.9]{Megginson1998}).

\begin{theorem}[Hahn--Banach]\label{Theo-HahnBanach}
Let $(\Omega,\mathcal{F},\mathbb{P})$ be a $\sigma$--generated probability space and $p\in(1,\infty)$.  
Given a countable dense subset $\{X_j:j\in\mathbb{N}\}$ of $L^{p}(\Omega)$, there exists a sequence of linear and continuous functionals $L_j:L^{p}(\Omega)\to\mathbb{R}$ such that
$$
L_j(X_j)=\|X_j\|_{p},\qquad \|L_j\|_{q}=1,
$$
and moreover
$$
\Vert X\Vert_{p}=\sup_{j\in\mathbb{N}}|L_j(X)|\qquad\forall\,X\in L^{p}(\Omega).
$$
\end{theorem}

In this way, the spaces $L^{p}(\Omega)$ constitute the natural setting in which random variables —now understood as elements of a function space— can be analyzed through the tools of functional analysis and measure theory.

Within this framework, the mapping $\chi:(0,1)\to L^{p}(0,1)$ from Problem \ref{Problem2} becomes the ideal model for exploring, in a rigorous way, \textit{what it means to be random in a Lebesgue space}.

\section{Constructing the Concept of a Random Variable in Lebesgue Spaces}

From now on, we shall consider $(\Omega,\mathcal{F},\mathbb{P})$ to be a complete and $\sigma$--generated probability space, and fix an exponent $p\in(1,\infty)$.

Recall first that, in the real case, the Borel $\sigma$--algebra $\mathcal{B}(\mathbb{R})$ is constructed from the open intervals $(a,b)$ with $a,b\in\mathbb{R}$ and $a<b$, although it can equivalently be obtained directly from the usual topology of $\mathbb{R}$. In other words, $\mathcal{B}(\mathbb{R})$ is the $\sigma$-algebra generated by the open sets determined by the metric induced by the absolute value.

Analogously, the $p$--norm in the Lebesgue space $L^{p}(\Omega)$ induces a natural topology and, with it, a corresponding Borel $\sigma$--algebra, which we denote by
$$
\mathcal{B}_{p}:=\mathcal{B}(L^{p}(\Omega)).
$$

Hence, the pair $(L^{p}(\Omega),\mathcal{B}_{p})$ becomes a measurable space. This structure allows us to extend the classical notion of measurability —and consequently, of a random variable— to the context of functions taking values in Banach spaces.

\begin{definition}\label{Def1_VALp}
A function $\xi:\Omega\to L^{p}(\Omega)$ is called an \textbf{$L^{p}(\Omega)$-valued random variable} if for every open set $\mathcal{O}$ on $L^{p}(\Omega)$ one has $\xi^{-1}(\mathcal{O})\in\mathcal{F}$.
\end{definition}

From this definition, the following result is immediate.

\begin{theorem}\label{Theo-Dis-Lp}
Let $\xi:\Omega\to L^{p}(\Omega)$ be a random variable.  
Then the function $\xi_{\sharp}\mathbb{P}:\mathcal{B}_{p}\to\mathbb{R}$ defined by
$$
\xi_{\sharp}\mathbb{P}(\mathcal{O}):=\mathbb{P}\big(\xi^{-1}(\mathcal{O})\big)
$$
defines a probability measure on $(L^{p}(\Omega),\mathcal{B}_{p})$.
\end{theorem}

\begin{proof}
The result follows directly from the properties of the measure $\mathbb{P}$ and from the stability of inverse images under countable unions and complements.
\end{proof}

Despite its formal elegance, Definition \ref{Def1_VALp} is not always practical for the study of concrete examples.  
This is mainly due to two reasons: first, the structure of the open subsets of $L^{p}(\Omega)$ can be difficult to describe explicitly; and second, even when they can be characterized in terms of open balls, determining their inverse images under $\xi$ may be technically cumbersome.

In the real case, there exist alternative characterizations of random variables. One of the most useful consists in describing them as limits of sequences of simple functions. Indeed, a function $X:\Omega\to\mathbb{R}$ is a random variable if and only if it is the almost sure pointwise limit of a sequence of simple random variables, that is, measurable functions taking finitely many values in $\mathbb{R}$ (see \cite[Result 1, Chapter 9]{JacodProtter2004}).

\begin{figure}[!ht]
\centering
\begin{tikzpicture}[xscale=2.5, yscale=2.5]
%\draw[fill=lime!20, dotted] (0,0)--plot[domain=0:4] (\x, {\x^(1/2)})--(4,0)--(0,0);
\draw[domain= 0:1.5, blue] plot(\x,{0.5*\x^2 -(1/3)*(\x)^3} );
\draw[domain= -1.5:0, blue] plot(\x,{-0.5*\x^2 -(1/6)*(\x)^6} );

\draw[-, gray] (-1.5,-0.78)--(1.522,-0.78);
\draw[-, gray] (-1.5,-0.78)--(-1.5,0.5);

\draw (1.5,-0.78) node[right]{$_{(\Omega,\mathcal{F})}$}; 
\draw (-1.5,0.5) node[left]{$_{(\mathbb{R},\mathcal{B}(\mathbb{R}))}$}; 

\draw[gray] (-1.5,-0.7) node{$_{-}$}; \draw[-, densely dotted,gray] (-1.5,-0.7)--(1.5,-0.7);
\draw[gray] (-1.5,-0.6) node{$_{-}$}; \draw[-, densely dotted,gray] (-1.5,-0.6)--(1.5,-0.6);
\draw[gray] (-1.5,-0.5) node{$_{-}$}; \draw[-, densely dotted,gray] (-1.5,-0.5)--(1.5,-0.5);
\draw[gray] (-1.5,-0.4) node{$_{-}$}; \draw[-, densely dotted,gray] (-1.5,-0.4)--(1.5,-0.4);
\draw[gray] (-1.5,-0.3) node{$_{-}$}; \draw[-, densely dotted,gray] (-1.5,-0.3)--(1.5,-0.3);
\draw[gray] (-1.5,-0.2) node{$_{-}$}; \draw[-, densely dotted,gray] (-1.5,-0.2)--(1.5,-0.2);
\draw[gray] (-1.5,-0.1) node{$_{-}$}; \draw[-, densely dotted,gray] (-1.5,-0.1)--(1.5,-0.1);
\draw[gray] (-1.5,0.0) node{$_{-}$}; \draw[-, densely dotted,gray] (-1.5,0.0)--(1.5,0.0);
\draw[gray] (-1.5,0.1) node{$_{-}$}; \draw[-, densely dotted,gray] (-1.5,0.1)--(1.5,0.1);
\draw[gray] (-1.5,0.2) node{$_{-}$}; \draw[-, densely dotted,gray] (-1.5,0.2)--(1.5,0.2);
\draw[gray] (-1.5,0.3) node{$_{-}$}; \draw[-, densely dotted,gray] (-1.5,0.3)--(1.5,0.3);

\draw[ultra thick] (-1.51,-0.78)--(-1.48,-0.78);  \draw[-, densely dotted] (-1.48,-0.7)--(-1.48,-0.78);
\draw[ultra thick] (-1.48,-0.7)--(-1.46,-0.7);  \draw[-, densely dotted] (-1.46,-0.6)--(-1.46,-0.7);
\draw[ultra thick] (-1.47,-0.6)--(-1.45,-0.6);  \draw[-, densely dotted] (-1.45,-0.5)--(-1.45,-0.6);
\draw[ultra thick] (-1.45,-0.5)--(-1.43,-0.5);  \draw[-, densely dotted] (-1.43,-0.4)--(-1.43,-0.5);
\draw[ultra thick] (-1.43,-0.4)--(-1.41,-0.4);  \draw[-, densely dotted] (-1.41,-0.3)--(-1.41,-0.4);
\draw[ultra thick] (-1.41,-0.3)--(-1.38,-0.3);  \draw[-, densely dotted] (-1.38,-0.2)--(-1.38,-0.3);
\draw[ultra thick] (-1.38,-0.2)--(-1.35,-0.2);  \draw[-, densely dotted] (-1.35,-0.1)--(-1.35,-0.2);
\draw[ultra thick] (-1.35,-0.1)--(-1.32,-0.1);  \draw[-, densely dotted] (-1.32,0)--(-1.32,-0.1);

\draw[ultra thick] (-1.32,0)--(-1.27,0);  \draw[-, densely dotted] (-1.27,0.1)--(-1.27,0);
\draw[ultra thick] (-0.45,0)--(0.57,0);  \draw[-, densely dotted] (-0.45,0.1)--(-0.45,0);
\draw[-, densely dotted] (0.57,0.1)--(0.57,0);
\draw[ultra thick] (1.32,0)--(1.5,0); \draw[-, densely dotted] (1.32,0.1)--(1.32,0);

\draw[ultra thick] (-1.27,0.1)--(-1.21,0.1);  \draw[-, densely dotted] (-1.21,0.2)--(-1.21,0.1);
\draw[ultra thick] (-0.65,0.1)--(-0.45,0.1);  \draw[-, densely dotted] (-0.65,0.2)--(-0.65,0.1);
\draw[ultra thick] (0.57,0.1)--(1.32,0.1);

\draw[ultra thick] (-1.21,0.2)--(-1.12,0.2);  \draw[-, densely dotted] (-1.12,0.3)--(-1.12,0.2);
\draw[ultra thick] (-0.85,0.2)--(-0.65,0.2);  \draw[-, densely dotted] (-0.85,0.3)--(-0.85,0.2);

\draw[ultra thick] (-0.85,0.3)--(-1.12,0.3);
\end{tikzpicture}
\end{figure}

This naturally suggests the following question: can Definition \ref{Def1_VALp} be characterized in terms of convergence of sequences of simple functions?

The answer was established by B.J.Pettis; cf. \cite{Pettis1938}, who formulated what is now known as the \textit{Pettis measurability theorem} (see \cite[Theorem 2, Chapter 2, Section 1]{Diestel}). To this end, he introduced several notions of measurability that directly connect the topological and linear structures of Banach spaces with measure theory.

\begin{definition}
Let $\xi:\Omega\to L^{p}(\Omega)$ be a function.
\begin{itemize}
    \item[(a)] $\xi$ is a \textbf{simple $L^{p}(\Omega)$--valued random variable} if there exist a finite partition $\{A_1,\ldots,A_N\}$ of $\Omega$, with $A_k\in\mathcal{F}$ for all $k=1,\ldots,N$, and functions $X_1,\ldots,X_N\in L^{p}(\Omega)$ such that
    $$
    \xi(\omega)=\sum_{k=1}^{N} X_k\,\mathbf{1}_{A_k}(\omega).
    $$

    \item[(b)] $\xi$ is \textbf{strongly measurable} if there exists a sequence of simple $L^{p}(\Omega)$--valued random variables $\xi_j:\Omega\to L^{p}(\Omega)$, $j\in\mathbb{N}$, such that $\xi_j\to\xi$ almost surely in $\Omega$.

    \item[(c)] $\xi$ is \textbf{weakly measurable} if for every linear and continuous functional $L:L^{p}(\Omega)\to\mathbb{R}$, the composition $L\circ\xi:\Omega\to\mathbb{R}$ is a real--valued random variable.
\end{itemize}
\end{definition}

These definitions lead to two fundamental questions: how are these notions related to each other, and what is their connection with the general definition of an $L^{p}(\Omega)$--valued random variable?

Let $\xi:\Omega\to L^{p}(\Omega)$ be a strongly measurable function. By definition, there exists a sequence of simple random variables $\xi_j:\Omega\to L^{p}(\Omega)$ such that
$$
\mathbb{P}\Big(\big\{\omega\in\Omega : \lim_{j\to\infty}\xi_j(\omega)=\xi(\omega)\big\}\Big)=1.
$$

Let $L:L^{p}(\Omega)\to\mathbb{R}$ be any linear and continuous functional.  

For each $j\in\mathbb{N}$, we may write
$$
\xi_j=\sum_{k=1}^{N_j} X_k^{(j)}\,\mathbf{1}_{A_k^{(j)}},
$$
where $X_1^{(j)},\ldots,X_{N_j}^{(j)}\in L^{p}(\Omega)$ and $\{A_1^{(j)},\ldots,A_{N_j}^{(j)}\}$ is a partition of $\Omega$ with $A_k^{(j)}\in\mathcal{F}$.

By the linearity of $L$, we have
$$
L(\xi_j(\omega))=\sum_{k=1}^{N_j}L\big(X_k^{(j)}\big)\,\mathbf{1}_{A_k^{(j)}}(\omega), \qquad \forall\,\omega\in\Omega,
$$
and therefore $L\circ\xi_j$ is a simple real--valued random variable.

Applying the continuity of $L$, it follows that
$$
\mathbb{P}\Big(\big\{\omega\in\Omega : \lim_{j\to\infty}L(\xi_j(\omega))=L(\xi(\omega))\big\}\Big)=1,
$$
which implies that $L\circ\xi$ is a real random variable. Consequently, every strongly measurable function is weakly measurable.

The converse question —whether every weakly measurable function is also strongly measurable— leads us to the central result of Pettis.

\begin{theorem}[Pettis Measurability Theorem]\label{Theo-M-Pettis}
Let $\xi:\Omega\to L^{p}(\Omega)$ be a function. The following statements are equivalent:
\begin{itemize}
    \item[(i)] $\xi$ is strongly measurable;
    \item[(ii)] $\xi$ is weakly measurable;
    \item[(iii)] $\xi$ is an $L^{p}(\Omega)$--valued random variable.
\end{itemize}
\end{theorem}

\begin{proof}
$(i)\Rightarrow(ii)$ was proved in the paragraph preceding the statement.

$(ii)\Rightarrow(iii)$: Since $L^{p}(\Omega)$ is separable (see Theorem \ref{Theo-Lp-Sep}), every open set $\mathcal{O}$ on $L^{p}(\Omega)$ can be written as a countable union of open balls in $L^{p}(\Omega)$. Therefore, it suffices to show that the inverse image under $\xi$ of any open ball in $L^{p}(\Omega)$ belongs to the $\sigma$-algebra $\mathcal{F}$ (see \cite[Theorem 8.1]{JacodProtter2004}).

Let $X\in L^{p}(\Omega)$ and $\varepsilon>0$ be arbitrary, and consider the open ball
$$
B_{p}(X,\varepsilon):=\bigl\{\,Y\in L^{p}(\Omega):\|X-Y\|_{p}<\varepsilon\,\bigr\}.
$$

Because $L^{p}(\Omega)$ is separable, fix a countable dense subset $\{X_{j}\in L^{p}(\Omega):j\in\mathbb{N}\}$. By the Hahn-Banach theorem (see Theorem \ref{Theo-HahnBanach}), there exists a sequence of linear and continuous functionals $L_{j}:L^{p}(\Omega)\to\mathbb{R}$ such that
$$
\|Z\|_{p}=\sup_{j\in\mathbb{N}}|L_{j}(Z)|\qquad\forall\,Z\in L^{p}(\Omega).
$$

Hence, for $Y\in L^{p}(\Omega)$, the condition $\|X-Y\|_{p}<\varepsilon$ is equivalent to the existence of some $k \in \mathbb{N}$ such that $|L_{j}(Y) - L_{j}(X)| < \varepsilon - \tfrac{1}{k}$ for all $j \in \mathbb{N}$.

Therefore,
$$
B_{p}(X,\varepsilon)=\bigcup_{k=1}^{\infty} \bigcap_{j =1}^{\infty} \left\{ Y \in L^{p}(\Omega)\,:\,|L_{j}(Y)-L_{j}(X)| <\varepsilon-\frac{1}{k}\right\}
$$
and consequently,
$$
\begin{aligned}
\xi^{-1}(B_{p}(X,\varepsilon))&=\xi^{-1}\left(\bigcup_{k=1}^{\infty} \bigcap_{j =1}^{\infty} \left\{ Y \in L^{p}(\Omega)\,:\,|L_{j}(Y)-L_{j}(X)| <\varepsilon-\frac{1}{k}\right\} \right)\\
&=\bigcup_{k=1}^{\infty}  \bigcap_{j =1}^{\infty} \xi^{-1}\left(\left\{ Y \in L^{p}(\Omega)\,:\,|L_{j}(Y)-L_{j}(X)| <\varepsilon-\frac{1}{k}\right\} \right)\\
&=\bigcup_{k=1}^{\infty}  \bigcap_{j =1}^{\infty} \left\{ \omega \in \Omega\,:\,|L_{j}(\xi(\omega))-L_{j}(X)| <\varepsilon-\frac{1}{k}\right\} \\
&=\bigcup_{k=1}^{\infty}  \bigcap_{j =1}^{\infty} (L_{j}\circ \xi)^{-1} \left(L_{j}(X)-\varepsilon+\frac{1}{k},L_{j}(X)+\varepsilon -\frac{1}{k}\right). \\
\end{aligned}
$$

Since $\xi$ is weakly measurable, each $L_{j}\circ\xi:\Omega\to\mathbb{R}$ is a real-valued random variable. Thus every set in the last line is in $\mathcal{F}$, hence $\xi^{-1}(B_{p}(X,\varepsilon))\in\mathcal{F}$. It follows that $\xi$ is an $L^{p}(\Omega)$--valued random variable.

$(iii)\Rightarrow(i)$: Let $\mathcal{D}=\{X_{j}\in L^{p}(\Omega):j\in\mathbb{N}\}$ be a countable dense subset of $L^{p}(\Omega)$. Then for each $X\in L^{p}(\Omega)$ and $k\in\mathbb{N}$ there exists $X_{j_{k}}\in\mathcal{D}$ such that
$$
\|X-X_{j_{k}}\|_{p}<\frac{1}{2^{k}},
$$
that is, $X\in B_{p}(X_{j_{k}},\tfrac{1}{2^{k}})$. In particular, for each fixed $k\in\mathbb{N}$,
$$
L^{p}(\Omega)=\bigcup_{j=1}^{\infty}B_{p}\left(X_{j},\frac{1}{2^{k}}\right).
$$

Fix $k\in\mathbb{N}$ and define, for $j\in\mathbb{N}$,
$$
E_{j,k}:=\bigl\{\omega\in\Omega:\ \xi(\omega)\in B_{p}(X_{j},2^{-k})\bigr\}=\xi^{-1}\!\bigl(B_{p}(X_{j},2^{-k})\bigr),
$$
which belongs to $\mathcal{F}$ since $\xi$ is an $L^{p}(\Omega)$--valued random variable. Define a disjoint sequence $\{A_{j,k} \in \mathcal{F}:j\in\mathbb{N}\}$ by
$$
A_{1,k}:=E_{1,k},\qquad
A_{j,k}:=E_{j,k}\smallsetminus\bigcup_{\ell=1}^{j-1}A_{\ell,k}\quad(j\ge2).
$$

Then $A_{j,k}\subset E_{j,k}$ and $\Omega=\bigcup_{j=1}^{\infty}A_{j,k}$. Hence, define $\xi_{k}:\Omega\to L^{p}(\Omega)$ by
$$
\xi_{k}(\omega):=\sum_{j=1}^{\infty}X_{j}\,\mathbf{1}_{A_{j,k}}(\omega).
$$

If $\omega\in A_{j,k}$, then $\xi_{k}(\omega)=X_{j}$ and, by the definition of $E_{j,k}$,
$$
\|\xi_{k}(\omega)-\xi(\omega)\|_{p}
=\|X_{j}-\xi(\omega)\|_{p}<\frac{1}{2^{k}}.
$$

Hence,
$$
\|\xi_{k}(\omega)-\xi(\omega)\|_{p}<\frac{1}{2^{k}}\qquad \forall\,\omega\in\Omega,
$$
and letting $k\to\infty$ yields
$$
\lim_{k\to\infty}\|\xi_{k}(\omega)-\xi(\omega)\|_{p}=0\qquad \forall\,\omega\in\Omega.
$$

Thus, $\xi_{k}\to\xi$ pointwise on $\Omega$.

Since $\{A_{j,k} \in \mathcal{F}\,:\, j \in \mathbb{N}\}$ is a partition of $\Omega$,
$$
\sum_{j=1}^{\infty}\mathbb{P}(A_{j,k})=1.
$$

By the Cauchy criterion for series, there exists $N_{k}\in\mathbb{N}$ such that
$$
\sum_{j=N_{k}+1}^{\infty}\mathbb{P}(A_{j,k})<\frac{1}{2^{k}}.
$$

Let $\widetilde{A}_{k}:=\bigcup_{j=N_{k}+1}^{\infty}A_{j,k}\in\mathcal{F}$ and define the truncation $\widetilde{\xi}_{k}:\Omega\to L^{p}(\Omega)$ by
$$
\widetilde{\xi}_{k}:=\sum_{j=1}^{N_{k}}X_{j}\,\mathbf{1}_{A_{j,k}}+X_{1}\,\mathbf{1}_{\widetilde{A}_{k}},
$$
which is clearly a simple $L^{p}(\Omega)$--valued random variable.  

If $\omega\in\Omega\smallsetminus\widetilde{A}_{k}$, then $\omega\notin A_{j,k}$ for all $j>N_{k}$, hence there is a unique $j_{\ast}\le N_{k}$ with $\omega\in A_{j_{*},k}$ and consequently $\xi_{k}(\omega)=X_{j_{\ast}}=\widetilde{\xi}_{k}(\omega)$. Therefore,
$$
\mathbb{P}\bigl(\{\omega\in\Omega:\ \widetilde{\xi}_{k}(\omega)\neq \xi_{k}(\omega)\}\bigr)\le \mathbb{P}(\widetilde{A}_{k})<\frac{1}{2^{k}}.
$$

Since $\sum_{k=1}^{\infty}\mathbb{P}(\widetilde{A}_{k})\le \sum_{k=1}^{\infty}\frac{1}{2^{k}}<\infty$, the Borel--Cantelli lemma (see \cite[Theorem 2.3.1]{Durrett2019}) implies that
$$
\mathbb{P}\!\left(\limsup_{k\to\infty}\widetilde{A}_{k}\right)=0.
$$

Let $A:=\Omega\smallsetminus\limsup_{k\to\infty}\widetilde{A}_{k}$. Then $\mathbb{P}(A)=1$ and, by the definition of the limsup,
$$
A=\Bigl\{\omega\in\Omega:\ \exists\,k_{0}\in\mathbb{N}\ \text{such that}\ \omega\in\Omega\smallsetminus \widetilde{A}_{k}\ \ \forall\,k\ge k_{0}\Bigr\}.
$$

If $\omega\in A$, there exists $k_{0} \in \mathbb{N}$ such that $\widetilde{\xi}_{k}(\omega)=\xi_{k}(\omega)$ for all $k\ge k_{0}$. Given $\varepsilon>0$, since $\xi_{k}\to\xi$ pointwise, there exists $k(\omega,\varepsilon) \in \mathbb{N}$ with
$$
\|\xi_{k}(\omega)-\xi(\omega)\|_{p}<\varepsilon\qquad \forall\,k\ge k(\omega,\varepsilon).
$$

Taking $k_{*}:=\max\{k_{0},k(\omega,\varepsilon)\} \in \mathbb{N}$, for all $k\ge k_{*}$ we have
$$
\widetilde{\xi}_{k}(\omega)=\xi_{k}(\omega)
\quad\text{and}\quad
\|\widetilde{\xi}_{k}(\omega)-\xi(\omega)\|_{p}<\varepsilon.
$$

Consequently,
$$
A\subset\bigl\{\omega\in\Omega:\ \widetilde{\xi}_{k}(\omega)\to \xi(\omega)\bigr\},
$$
and $1\leq \mathbb{P}(A)\leq \mathbb{P}\left(\bigl\{\omega\in\Omega:\ \widetilde{\xi}_{k}(\omega)\to \xi(\omega)\bigr\}\right)\leq 1$ which shows that $\widetilde{\xi}_{k}\to\xi$ almost surely on $\Omega$. Therefore, $\xi$ is strongly measurable.
\end{proof}

The preceding theorem completes the characterization of an $L^{p}(\Omega)$--valued random variable, showing that Pettis measurability —in its strong, weak, and Borel forms— coincides in this setting. Thus, the initial definition, apparently abstract, translates into properties formulated in familiar terms of functional analysis, such as density, duality, and the continuity of linear functionals. This perspective not only clarifies the topological and probabilistic structure of vector-valued random variables but also provides more flexible tools for their study, especially within Banach spaces; cf. \cite{Diestel,Pettis1938}.

The characterization of random variables through the Pettis measurability theorem naturally suggests a further direction: once it is understood what it means for a function to take random values in a Lebesgue space, the next question is how to define its expected value, that is, how to extend the notion of expectation to functions that are no longer real-valued but elements of $L^{p}(\Omega)$. In the following section, we address precisely this question, showing that the notion of $L^{p}(\Omega)$--valued expectation can be coherently formulated via the Bochner integral, preserving the functional and geometric interpretation of the classical concept.

\section{The Concept of Expectation for Random Variables in Lebesgue Spaces}

In classical probability theory, the expectation of a real random variable $X:\Omega\to\mathbb{R}$ is defined through the Lebesgue integral with respect to the probability measure. This integral not only assigns a number summarizing the average behavior of $X$, but also embodies a deeper geometric and analytic idea: the expectation is the probabilistic center of mass of the values taken by the variable.

From this point of view, computing the expectation amounts to finding the point that balances the ``weight'' of the distribution of $X$ along the real line. This interpretation, beyond its algebraic role as a linear functional, endows the expectation with a structural meaning that transcends the mere operation of integration.

However, when random variables take values in function spaces such as $L^{p}(\Omega)$, the situation changes radically. It is no longer possible to speak of averaging real coordinates, and the notion of a ``center'' must be reinterpreted in a setting where the values of the random variable are themselves functions. A natural question therefore arises: \textit{how should one define the expectation of an $L^{p}(\Omega)$--valued random variable?}

The path toward answering this question was opened in the early twentieth century, when mathematicians sought an extension of the Lebesgue integral to vector--valued functions. This generalization, formulated in the context of Banach spaces, became known as the \textit{Bochner integral}; cf. \cite{Diestel,Pettis1938}.

The essential idea is to preserve the fundamental properties of the real integral —linearity, continuity, and monotone and dominated convergence— but reformulated in terms of the norm of the space in which the function takes its values. Instead of requiring pointwise integrability, one demands that the norm be integrable, thereby quantifying the average size of the random vector and ensuring that the integral defines a well-defined element of the same vector space.

Historically, this extension was developed for functions taking values in separable Banach spaces, where the metric structure and duality guarantee the existence of a coherent notion of measurability and vector integration.

The original construction, due to Bochner and elaborated in subsequent works by Pettis and Dunford; cf. \cite{Pettis1938}, is technically delicate. However, in the specific context of Lebesgue spaces $L^{p}(\Omega)$, it can be presented in a more direct and practical form, exploiting the Riesz representation theorem.

In this section, we develop a simplified construction of the concept of expectation for random variables taking values in $L^{p}(\Omega)$ spaces, avoiding the full machinery of the Bochner integral while retaining its essential analytic meaning.

Let $\xi:\Omega \to L^{p}(\Omega)$ be a random variable. To analyze the size of the values taken by $\xi$, consider the norm function $\|\cdot\|_{p}:L^{p}(\Omega)\to\mathbb{R}$. Since the norm is continuous, it is measurable, and consequently the composition $\|\cdot\|_{p}\circ\xi:\Omega\to\mathbb{R}$ defines a nonnegative real-valued random variable.  
This allows us to quantify the magnitude of the values of $\xi$ and to establish integrability conditions in $L^{p}(\Omega)$.

Suppose therefore that $\xi:\Omega \to L^{p}(\Omega)$ satisfies
$$
\int_{\Omega}\|\xi(\omega)\|_{p}\,d\mathbb{P}(\omega)<\infty.
$$

For $p,q\in(1,\infty)$ with $\tfrac{1}{p}+\tfrac{1}{q}=1$, Hölder’s inequality (see Theorem \ref{HolderInequality}) allows us to define the function 
$\beta:L^{p}(\Omega)\times L^{q}(\Omega)\to\mathbb{R}$ by
$$
\beta(X,Y):=\int_{\Omega}X(\upsilon)Y(\upsilon)\,d\mathbb{P}(\upsilon).
$$

By the linearity of the Lebesgue integral, we have
$$
\begin{aligned}
\beta\!\left(aX_{1}+X_{2},\,Y_{1}+bY_{2}\right)
&=\int_{\Omega}(aX_{1}(\upsilon)+X_{2}(\upsilon))(Y_{1}(\upsilon)+bY_{2}(\upsilon))\,d\mathbb{P}(\upsilon)\\
&=\int_{\Omega}\!\big(aX_{1}(\upsilon)Y_{1}(\upsilon)+abX_{1}(\upsilon)Y_{2}(\upsilon)\big)\,d\mathbb{P}(\upsilon)\\
&\quad+\int_{\Omega}\!\big(X_{2}(\upsilon)Y_{1}(\upsilon)+bX_{2}(\upsilon)Y_{2}(\upsilon)\big)\,d\mathbb{P}(\upsilon)\\
&=a\,\beta(X_{1},Y_{1})+ab\,\beta(X_{1},Y_{2})+\beta(X_{2},Y_{1})+b\,\beta(X_{2},Y_{2})
\end{aligned}
$$
for any $X_{1},X_{2}\in L^{p}(\Omega)$, $Y_{1},Y_{2}\in L^{q}(\Omega)$, and $a,b\in\mathbb{R}$.  

Moreover, Hölder’s inequality ensures that
$$
|\beta(X,Y)|\le\int_{\Omega}|X(\upsilon)Y(\upsilon)|\,d\mathbb{P}(\upsilon)
\le\|X\|_{p}\|Y\|_{q},
$$
for all $X\in L^{p}(\Omega)$ and $Y\in L^{q}(\Omega)$.

In particular, for each fixed $Y\in L^{q}(\Omega)$, the mapping $\beta(\,\cdot\,,Y):L^{p}(\Omega)\to\mathbb{R}$ defines a linear and continuous functional.

Since $\xi:\Omega\to L^{p}(\Omega)$ is a random variable, the Pettis measurability theorem (see Theorem~\ref{Theo-M-Pettis}) ensures that the function $\beta(\xi,Y):\Omega\to\mathbb{R}$ defined by
$$
\beta(\xi,Y)(\omega)=\beta\big(\xi(\omega),Y\big)
$$
is a real-valued random variable.  

Applying Hölder’s inequality once again, we obtain
$$
0\le|\beta(\xi(\omega),Y)|\le\|\xi(\omega)\|_{p}\|Y\|_{q}
\qquad\forall\,\omega\in\Omega,
$$
and by the monotonicity of the Lebesgue integral,
$$
0\le\int_{\Omega}|\beta(\xi(\omega),Y)|\,d\mathbb{P}(\omega)
\le\int_{\Omega}\|\xi(\omega)\|_{p}\|Y\|_{q}\,d\mathbb{P}(\omega)
=\|Y\|_{q}\int_{\Omega}\|\xi(\omega)\|_{p}\,d\mathbb{P}(\omega)<\infty,
$$
so that the integral on the left-hand side is well defined.

This allows us to define the functional $\Phi:L^{q}(\Omega)\to\mathbb{R}$ by
$$
\Phi(Y):=\int_{\Omega}\beta\big(\xi(\omega),Y\big)\,d\mathbb{P}(\omega).
$$

The mapping $\Phi$ is clearly linear. Using Hölder’s inequality (see Theorem \ref{HolderInequality}) we have
$$
|\Phi(Y)|\le\|Y\|_{q}\int_{\Omega}\|\xi(\omega)\|_{p}\,d\mathbb{P}(\omega)
\qquad\forall\,Y\in L^{q}(\Omega),
$$
which shows that $\Phi$ is continuous on $L^{q}(\Omega)$.  

By the Riesz representation theorem (see Theorem \ref{The-Representation-Riesz}), there exists a unique random variable $\mathcal{I}\in L^{p}(\Omega)$ such that
$$
\Phi(Y)=\int_{\Omega}\mathcal{I}(\upsilon)Y(\upsilon)\,d\mathbb{P}(\upsilon)
\qquad\forall\,Y\in L^{q}(\Omega).
$$

We call this element $\mathcal{I}$ the \textbf{Bochner integral} of the random variable $\xi:\Omega\to L^{p}(\Omega)$, and, for convenience, we denote it by $\mathcal{I}=\mathbb{E}(\xi)$.

\begin{definition}
Let $\xi:\Omega\to L^{p}(\Omega)$ be a random variable.  
We say that $\xi$ is \textbf{Bochner--integrable} if
$$
\int_{\Omega}\|\xi(\omega)\|_{p}\,d\mathbb{P}(\omega)<\infty,
$$
and we define its \textbf{expectation} (Bochner integral) as the unique random variable $\mathbb{E}(\xi)\in L^{p}(\Omega)$ satisfying
$$
\int_{\Omega}\beta\big(\xi(\omega),Y\big)\,d\mathbb{P}(\omega)
=\int_{\Omega}\mathbb{E}(\xi)(\upsilon)\,Y(\upsilon)\,d\mathbb{P}(\upsilon)
\qquad\forall\,Y\in L^{q}(\Omega).
$$
\end{definition}

This construction shows that the expectation $\mathbb{E}(\xi)$ is the unique element of $L^{p}(\Omega)$ reproducing, for every functional in the dual space $L^{q}(\Omega)$, the average action of the integration functional induced by $\xi$. In other words, the Bochner integral provides the Riesz representation of the functional $\Phi$, thereby extending the classical notion of expectation to the setting of $L^{p}(\Omega)$--valued random variables in a natural and geometrically meaningful way (see \cite[Theorems 2 and 6, Chapter 2, Section 2]{Diestel}).

\section{Solution to the Initial Problem}

Having established, from a functional perspective, the concept of a random variable taking values in the Lebesgue spaces $L^{p}$, and having developed a coherent notion of expectation through the Bochner integral, we are now equipped with the necessary tools to address the problem posed at the beginning of this work (see Problem \ref{Problem2}).  

In particular, the framework constructed above allows us to interpret, in a rigorous way, the mapping $\chi:(0,1)\to L^{p}(0,1)$,
$$
\chi(\omega)=\mathbf{1}_{(0,\omega)},
$$
as a concrete example of an $L^{p}(0,1)$-valued random variable, and to compute its expectation within the very space of functions.

Let $\varepsilon\in(0,1)$ and $p\in(1,\infty)$.  
The essential idea is the following: once it has been established that the function $\chi:(0,1)\to L^{p}(0,1)$ is indeed a random variable, we will be able to determine the probability of the set
$$
\left\{\omega\in(0,1): \int_{0}^{1}\mathbf{1}_{(0,\omega)}(t)\,dt < \varepsilon^{p}\right\},
$$
by means of the distribution induced by $\chi$ (see Theorem \ref{Theo-Dis-Lp}). In particular, we shall look for a Borel subset $\mathcal{O}$ of $L^{p}(0,1)$ such that
$$
\chi^{-1}(\mathcal{O})
=\left\{\omega\in(0,1): \int_{0}^{1}\mathbf{1}_{(0,\omega)}(t)\,dt < \varepsilon^{p}\right\},
$$
which will allow us to apply the Lebesgue measure directly.

\begin{theorem}\label{Teo-SolucionProblemaInicial}
The function $\chi:(0,1)\to L^{p}(0,1)$ defined by
$$
\chi(\omega)=\mathbf{1}_{(0,\omega)},
$$
is an $L^{p}(0,1)$--valued random variable.
\end{theorem}

Since $((0,1),\mathcal{L}(0,1),\lambda)$ is a complete and $\sigma$--generated probability space, the Lebesgue space $L^{p}(0,1)$ is separable. Consequently, by the Pettis measurability theorem, Theorem \ref{Teo-SolucionProblemaInicial} can be proved in three equivalent ways:
\begin{itemize}
    \item[(1)] by showing that $\chi^{-1}(\mathcal{O})\in\mathcal{L}(0,1)$ for every open set $\mathcal{O}\subset L^{p}(0,1)$;
    \item[(2)] by verifying that $\chi$ is weakly measurable;
    \item[(3)] by showing that $\chi$ is strongly measurable.
\end{itemize}

In this simple example we shall consider all three approaches, both to illustrate the concepts introduced and to let the reader appreciate which provides the most direct path.

\begin{proof}[First proof of Theorem~\ref{Teo-SolucionProblemaInicial}]
We show directly that the function $\chi:(0,1)\to L^{p}(0,1)$ is continuous and, consequently, that it is an $L^{p}(0,1)$--valued random variable.

Let $\omega_{0}\in(0,1)$ and $\varepsilon>0$. Set $\delta:=\varepsilon^{p}>0$. Then, for every $\omega\in(0,1)$ such that $0<|\omega-\omega_{0}|<\delta$, we have
$$
\begin{aligned}
\|\chi(\omega_{0})-\chi(\omega)\|_{p}^{p}
&=\int_{(0,1)}|\mathbf{1}_{(0,\omega_{0})}(t)-\mathbf{1}_{(0,\omega)}(t)|^{p}\,d\lambda(t)\\[4pt]
&=\int_{(0,1)}|\mathbf{1}_{(0,\omega_{0})\triangle(0,\omega)}(t)|^{p}\,d\lambda(t)\\[4pt]
&=\lambda\big((0,\omega_{0})\triangle(0,\omega)\big)\\[4pt]
&=|\omega_{0}-\omega|\\[4pt]
&<\varepsilon^{p}.
\end{aligned}
$$

Therefore, $\chi$ is continuous on $(0,1)$.
\end{proof}

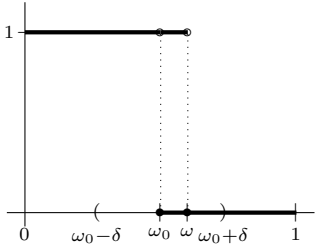
\begin{figure}[!ht]
\centering
\begin{tikzpicture}[xscale=1.2,yscale=0.8]
	\draw[-] (0,0) -- (3.2,0); \draw [-] (0,0) -- (0,3.5); 
	\draw[-] (-0.2,0)--(0,0); \draw[-] (0,0)--(0,-0.23);

\draw (0,3) node{$_{-}$}; \draw (0,3) node[left]{$_{1}$};
\draw (1.5,0) node{$_{|}$}; \draw (1.5,-0.1) node[below]{$_{\omega_{0}}$};
\draw (0.8,-0.1) node[below]{$_{\omega_{0}-\delta}$};
\draw (2.2,-0.1) node[below]{$_{\omega_{0}+\delta}$};
\draw (0.8,0) node{$_{(}$};
\draw (2.2,0) node{$_{)}$};
\draw (3,0) node{$_{|}$}; \draw (3,-0.1) node[below]{$_{1}$};
\draw[ultra thick] (1.5,0)--(3,0); \draw (1.5,0) node{$_{\bullet}$}; \draw (1.5,3) node{$_{\circ}$};
\draw[ultra thick] (0,3)--(1.5,3); 
\draw (0,-0.1) node[below]{$_{0}$};
\draw [dotted] (1.5,0)--(1.5,3);
%\draw [dotted] (0.8,0)--(0.8,3);
%\draw [dotted] (2.2,0)--(2.2,3);
\draw (1.8,-0.1) node[below]{$_{\omega}$};
\draw (1.8,0) node{$_{|}$};
\draw[ultra thick] (1.5,3)--(1.8,3); 
\draw [dotted] (1.8,0)--(1.8,3);
\draw (1.8,0) node{$_{\bullet}$}; \draw (1.8,3) node{$_{\circ}$};
\end{tikzpicture}
%\caption{Gráfica de la función $f(x):=1_{(0,x)}$}
\end{figure}

\begin{proof}[Second proof of Theorem \ref{Teo-SolucionProblemaInicial}]
We shall prove that the function $\chi:(0,1)\to L^{p}(0,1)$ is weakly measurable.

Let $L:L^{p}(0,1)\to\mathbb{R}$ be a linear and continuous functional.  
By the Riesz representation theorem (see Theorem~\ref{The-Representation-Riesz}), there exists a unique function $Y\in L^{q}(0,1)$ such that
$$
L(\chi(\omega))=\int_{(0,1)}Y(t)\,\mathbf{1}_{(0,\omega)}(t)\,d\lambda(t),
\qquad \forall\,\omega\in(0,1).
$$

For any $\omega_{0},\omega\in(0,1)$ we have
$$
\begin{aligned}
|L(\chi(\omega_{0}))-L(\chi(\omega))|
&=|L(\chi(\omega_{0})-\chi(\omega))|\\[4pt]
&\leq \int_{(0,1)}|Y(t)|\,|\mathbf{1}_{(0,\omega_{0})}(t)-\mathbf{1}_{(0,\omega)}(t)|\,d\lambda(t)\\[4pt]
&\leq \|Y\|_{q}\,\|\mathbf{1}_{(0,\omega_{0})\triangle(0,\omega)}\|_{p}\\[4pt]
&\leq \|Y\|_{q}\,|\omega_{0}-\omega|^{1/p}.
\end{aligned}
$$

Hence, given $\varepsilon>0$, by defining $\delta^{1/p}:=\frac{\varepsilon}{\|Y\|_{q}+1}>0$, we obtain
$$
|L(\chi(\omega_{0}))-L(\chi(\omega))|<\varepsilon 
\qquad \text{whenever}\quad 0<|\omega-\omega_{0}|<\delta.
$$

Consequently, the composition $L\circ\chi:(0,1)\to\mathbb{R}$ is continuous, and thus a real-valued random variable.  

Since the functional $L$ is arbitrary, we conclude that $\chi$ is weakly measurable.
\end{proof}

\begin{proof}[Third proof of Theorem \ref{Teo-SolucionProblemaInicial}]
We construct a sequence $(\chi_{k})$ of simple $L^{p}(0,1)$--valued random variables that converges almost surely to the function $\chi$ on $(0,1)$.

For each $k\in\mathbb{N}$, partition the interval $(0,1]$ into the subintervals
$$
I_{j,k}:=\Big(\tfrac{j-1}{2^{k}},\,\tfrac{j}{2^{k}}\Big], 
\qquad j=1,\dots,2^{k}.
$$

Define $\chi_{j,k}:=\mathbf{1}_{(0,\,j/2^{k})}\in L^{p}(0,1)$ and the function $\chi_{k}:(0,1)\to L^{p}(0,1)$ by
$$
\chi_{k}(\omega):=\sum_{j=1}^{2^{k}}\chi_{j,k}\,\mathbf{1}_{I_{j,k}}(\omega),
$$
which, by construction, is a simple $L^{p}(0,1)$--valued random variable.

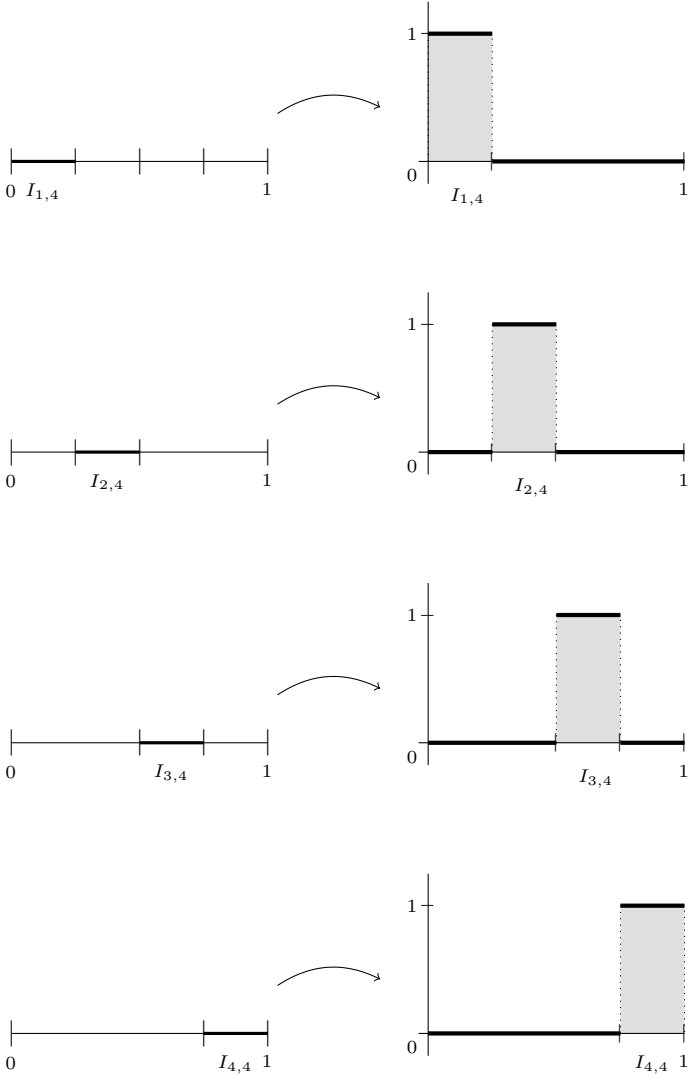
\begin{figure}[!ht]
    \centering
	\begin{tikzpicture}[xscale=0.85,yscale=0.85]
	\draw[-] (0,0)--(4,0);
	\draw (0,0) node{$|$}; \draw (0,-0.22) node[below]{$_{0}$};
	\draw (4,0) node{$|$}; \draw (4,-0.2) node[below]{$_{1}$};
    
	\draw (1,0) node{$|$}; \draw (0.5,-0.2) node[below]{$_{I_{1,4}}$};
    \draw (2,0) node{$|$}; %\draw (2,-0.2) node[below]{$_{I_{2,4}}$};
    \draw (3,0) node{$|$}; %\draw (3,-0.2) node[below]{$_{I_{3,4}}$};

    \draw[-,very thick] (0,0)--(1,0);
	%\draw[->] (0,-2)--(0,2);
	%\draw[-, densely dashed] (0,1)--(2.35,1);
%\draw[-, densely dashed] (0,-1)--(2.35,-1);
%\draw[domain= 0:2.35, thick, blue] plot(\x,{\x^(1/2)} );

    \draw[<-] (5.75,0.85) to[bend right] (4.15,0.75);
    %\draw(4.5,1.25)node[above]{$_{\chi}$};
    
    \draw [fill=gray!25,  dotted]  (6.5,0)--(6.5,2)--(7.5,2)--(7.5,0);
    
	\draw[-] (6.35,0)--(10.5,0);
	\draw[-] (6.5,-0.35)--(6.5,2.5);
\draw (6.5,2) node{$_{-}$}; \draw (6.5,2) node[left] {$_{1}$}; \draw (6.5,-0.22) node[left] {$_{0}$};
\draw (7.5,0) node{$_{|}$}; \draw (7.125,-0.27) node[below] {$_{I_{1,4}}$};
\draw (10.5,0) node{$_{|}$}; \draw (10.5,-0.2) node[below] {$_{1}$};
\draw [ultra thick,black] (7.5,0)--(10.5,0);
\draw [ultra thick,black] (6.5,2)--(7.5,2);
\end{tikzpicture}
%\caption{Representación gráfica}
\end{figure}

\begin{figure}[!ht]
    \centering
	\begin{tikzpicture}[xscale=0.85,yscale=0.85]
	\draw[-] (0,0)--(4,0);
	\draw (0,0) node{$|$}; \draw (0,-0.22) node[below]{$_{0}$};
	\draw (4,0) node{$|$}; \draw (4,-0.2) node[below]{$_{1}$};
    
	\draw (1,0) node{$|$}; %\draw (0.5,-0.2) node[below]{$_{I_{1,4}}$};
    \draw (2,0) node{$|$}; \draw (1.5,-0.2) node[below]{$_{I_{2,4}}$};
    %\draw (3,0) node{$|$}; \draw (3,-0.2) node[below]{$_{I_{3,4}}$};

    \draw[-,very thick] (1,0)--(2,0);
	%\draw[->] (0,-2)--(0,2);
	%\draw[-, densely dashed] (0,1)--(2.35,1);
%\draw[-, densely dashed] (0,-1)--(2.35,-1);
%\draw[domain= 0:2.35, thick, blue] plot(\x,{\x^(1/2)} );

    \draw[<-] (5.75,0.85) to[bend right] (4.15,0.75);
    %\draw(4.5,1.25)node[above]{$_{\chi}$};
    
    \draw [fill=gray!25,  dotted]  (7.5,0)--(7.5,2)--(8.5,2)--(8.5,0);
    
	\draw[-] (6.35,0)--(10.5,0);
	\draw[-] (6.5,-0.35)--(6.5,2.5);
\draw (6.5,2) node{$_{-}$}; \draw (6.5,2) node[left] {$_{1}$}; \draw (6.5,-0.22) node[left] {$_{0}$};
\draw (7.5,0) node{$_{|}$}; \draw (8.125,-0.27) node[below] {$_{I_{2,4}}$};
\draw (10.5,0) node{$_{|}$}; \draw (10.5,-0.2) node[below] {$_{1}$};

\draw (8.5,0) node{$_{|}$};

\draw [ultra thick,black] (8.5,0)--(10.5,0);
\draw [ultra thick,black] (6.5,0)--(7.5,0);

\draw [ultra thick,black] (7.5,2)--(8.5,2);

\end{tikzpicture}
%\caption{Representación gráfica}
\end{figure}

\begin{figure}[!ht]
    \centering
	\begin{tikzpicture}[xscale=0.85,yscale=0.85]
	\draw[-] (0,0)--(4,0);
	\draw (0,0) node{$|$}; \draw (0,-0.22) node[below]{$_{0}$};
	\draw (4,0) node{$|$}; \draw (4,-0.2) node[below]{$_{1}$};
    
	\draw (2,0) node{$|$}; %\draw (0.5,-0.2) node[below]{$_{I_{1,4}}$};
    \draw (3,0) node{$|$}; \draw (2.5,-0.2) node[below]{$_{I_{3,4}}$};
    %\draw (3,0) node{$|$}; \draw (3,-0.2) node[below]{$_{I_{3,4}}$};

    \draw[-,very thick] (2,0)--(3,0);
	%\draw[->] (0,-2)--(0,2);
	%\draw[-, densely dashed] (0,1)--(2.35,1);
%\draw[-, densely dashed] (0,-1)--(2.35,-1);
%\draw[domain= 0:2.35, thick, blue] plot(\x,{\x^(1/2)} );

    \draw[<-] (5.75,0.85) to[bend right] (4.15,0.75);
    %\draw(4.5,1.25)node[above]{$_{\chi}$};
    
    \draw [fill=gray!25,  dotted]  (8.5,0)--(8.5,2)--(9.5,2)--(9.5,0);
    
	\draw[-] (6.35,0)--(10.5,0);
	\draw[-] (6.5,-0.35)--(6.5,2.5);
\draw (6.5,2) node{$_{-}$}; \draw (6.5,2) node[left] {$_{1}$}; \draw (6.5,-0.22) node[left] {$_{0}$};
\draw (8.5,0) node{$_{|}$}; 
\draw (9.5,0) node{$_{|}$};

\draw (9.125,-0.27) node[below] {$_{I_{3,4}}$};
\draw (10.5,0) node{$_{|}$}; \draw (10.5,-0.2) node[below] {$_{1}$};

\draw [ultra thick,black] (9.5,0)--(10.5,0);
\draw [ultra thick,black] (6.5,0)--(8.5,0);

\draw [ultra thick,black] (8.5,2)--(9.5,2);

\end{tikzpicture}
%\caption{Representación gráfica}
\end{figure}

\begin{figure}[!ht]
    \centering
	\begin{tikzpicture}[xscale=0.85,yscale=0.85]
	\draw[-] (0,0)--(4,0);
	\draw (0,0) node{$|$}; \draw (0,-0.22) node[below]{$_{0}$};
	\draw (4,0) node{$|$}; \draw (4,-0.2) node[below]{$_{1}$};
    
	%\draw (2,0) node{$|$}; %\draw (0.5,-0.2) node[below]{$_{I_{1,4}}$};
    \draw (3,0) node{$|$}; \draw (3.5,-0.2) node[below]{$_{I_{4,4}}$};
    %\draw (3,0) node{$|$}; \draw (3,-0.2) node[below]{$_{I_{3,4}}$};

    \draw[-,very thick] (3,0)--(4,0);
	%\draw[->] (0,-2)--(0,2);
	%\draw[-, densely dashed] (0,1)--(2.35,1);
%\draw[-, densely dashed] (0,-1)--(2.35,-1);
%\draw[domain= 0:2.35, thick, blue] plot(\x,{\x^(1/2)} );

    \draw[<-] (5.75,0.85) to[bend right] (4.15,0.75);
    %\draw(4.5,1.25)node[above]{$_{\chi}$};
    
    \draw [fill=gray!25,  dotted]  (9.5,0)--(9.5,2)--(10.5,2)--(10.5,0);
    
	\draw[-] (6.35,0)--(10.5,0);
	\draw[-] (6.5,-0.35)--(6.5,2.5);
\draw (6.5,2) node{$_{-}$}; \draw (6.5,2) node[left] {$_{1}$}; \draw (6.5,-0.22) node[left] {$_{0}$};

\draw (9.5,0) node{$_{|}$}; 
%\draw (9.5,0) node{$_{|}$};

\draw (10,-0.2) node[below] {$_{I_{4,4}}$};
\draw (10.5,0) node{$_{|}$}; \draw (10.5,-0.2) node[below] {$_{1}$};

\draw [ultra thick,black] (9.5,2)--(10.5,2);

\draw [ultra thick,black] (6.5,0)--(9.5,0);

%\draw [ultra thick,black] (9.5,2)--(10.5,2);

\end{tikzpicture}
%\caption{Representación gráfica}
\end{figure}

If $\omega\in I_{j,k}$, then $\chi(\omega)=\mathbf{1}_{(0,\omega)}$ and $\chi_{k}(\omega)=\mathbf{1}_{(0,\,j/2^{k})}$.  
Hence,
$$
\|\chi_{k}(\omega)-\chi(\omega)\|_{p}
=\big\|\mathbf{1}_{(0,\,j/2^{k})}-\mathbf{1}_{(0,\omega)}\big\|_{p}
=\left|\tfrac{j}{2^{k}}-\omega\right|^{1/p}
\le \frac{1}{2^{k/p}}.
$$

Therefore,
$$
\lim_{k\to\infty}\|\chi_{k}(\omega)-\chi(\omega)\|_{p}=0,
\qquad \forall\,\omega\in(0,1).
$$

Consequently, the sequence $(\chi_{k})$ converges pointwise in the $L^{p}(0,1)$--norm to $\chi$. It follows that $\chi$ is strongly measurable.
\end{proof}

Once it has been established that $\chi:(0,1)\to L^{p}(0,1)$ is a random variable, we are ready to answer Problem \ref{Problem2}.

Let $B_{p}(0,\varepsilon)$ denote the open ball in $L^{p}(0,1)$ centered at the constant function $0$ with radius $\varepsilon\in(0,1)$. Then,
$$
\begin{aligned}
\chi^{-1}(B_{p}(0,\varepsilon))
&=\{\omega\in(0,1):\|\chi(\omega)\|_{p}<\varepsilon\}\\[4pt]
&=\left\{\omega\in(0,1):\left(\int_{(0,1)}\mathbf{1}_{(0,\omega)}(t)\,d\lambda(t)\right)^{1/p}<\varepsilon\right\}\\[4pt]
&=\left\{\omega\in(0,1):\int_{(0,1)}\mathbf{1}_{(0,\omega)}(t)\,d\lambda(t)<\varepsilon^{p}\right\}\\[4pt]
&=\{\omega\in(0,1):\omega<\varepsilon^{p}\}\\[4pt]
&=(0,\varepsilon^{p}).
\end{aligned}
$$

Therefore,
$$
\lambda\!\left(\chi^{-1}(B_{p}(0,\varepsilon))\right)
=\lambda((0,\varepsilon^{p}))=\varepsilon^{p}.
$$

\begin{figure}[!ht]
    \centering
	\begin{tikzpicture}[xscale=0.85,yscale=0.85]
    %\draw(4.5,1.25)node[above]{$_{\chi}$};
    \draw [fill=gray!25,  dotted]  (6.5,0)--(6.5,2)--(8,2)--(8,0);
	\draw[-] (6.35,0)--(9.35,0);
	\draw[-] (6.5,-0.35)--(6.5,2.5);
\draw (6.5,2) node{$_{-}$}; \draw (6.5,2) node[left] {$_{1}$}; \draw (6.5,-0.22) node[left] {$_{0}$};
\draw (8,0) node{$_{|}$}; \draw (8,-0.27) node[below] {$_{\omega}$};
\draw (8.45,0) node{$_{|}$}; \draw (8.45,-0.22) node[below] {$_{\varepsilon^{p}}$};
\draw (9,0) node{$_{|}$}; \draw (9,-0.2) node[below] {$_{1}$};
\draw [ultra thick,black] (8,0)--(9,0);
\draw [ultra thick,black] (6.5,2)--(8,2);

\end{tikzpicture}
%\caption{Representación gráfica}
\end{figure}

Before concluding, let us compute the expectation of the random variable $\chi:(0,1)\to L^{p}(0,1)$.
We observe that
$$
\int_{(0,1)}\|\chi(\omega)\|_{p}\,d\lambda(\omega)
=\int_{(0,1)}\omega^{p}\,d\lambda(\omega)
=\frac{1}{p+1}<\infty,
$$
so $\chi$ is Bochner integrable.

Fix an arbitrary $Y\in L^{q}(0,1)$. Applying Fubini’s theorem (see \cite[Theorem~1.7.2]{Durrett2019}), we obtain
$$
\begin{aligned}
\int_{(0,1)}\beta(\chi(\omega),Y)\,d\lambda(\omega)
&=\int_{(0,1)}\left(\int_{(0,1)}\mathbf{1}_{(0,\omega)}(t)\,Y(t)\,d\lambda(t)\right)\,d\lambda(\omega)\\[4pt]
&=\int_{(0,1)}\left(\int_{(0,1)}\mathbf{1}_{(t,1)}(\omega)\,Y(t)\,d\lambda(\omega)\right)\,d\lambda(t)\\[4pt]
&=\int_{(0,1)}Y(t)\left(\int_{(0,1)}\mathbf{1}_{(t,1)}(\omega)\,d\lambda(\omega)\right)\,d\lambda(t)\\[4pt]
&=\int_{(0,1)}(1-t)\,Y(t)\,d\lambda(t).
\end{aligned}
$$

By the Riesz representation theorem, the expectation (Bochner integral) of the random variable $\chi$ is the function $\mathbb{E}(\chi):(0,1)\to\mathbb{R}$ given by
$$
\mathbb{E}(\chi)(t)=1-t,
$$
which belongs to $L^{p}(0,1)$.

\begin{figure}[!ht]
    \centering
	\begin{tikzpicture}[xscale=0.85,yscale=0.85]
	\draw[-] (0,0)--(2.5,0);
	\draw (0,0) node{$|$}; \draw (0,-0.22) node[below]{$_{0}$};
	\draw (2.5,0) node{$|$}; \draw (2.5,-0.2) node[below]{$_{1}$};
	\draw (1.35,0) node{$|$}; \draw (1.35,-0.2) node[below]{$_{\omega}$};
	%\draw[->] (0,-2)--(0,2);
	%\draw[-, densely dashed] (0,1)--(2.35,1);
%\draw[-, densely dashed] (0,-1)--(2.35,-1);
%\draw[domain= 0:2.35, thick, blue] plot(\x,{\x^(1/2)} );

    \draw[<-] (5.75,0.85) to[bend right] (3.15,0.75);
    \draw(4.5,1.25)node[above]{$_{\chi}$};
    \draw [fill=gray!25,  dotted]  (6.5,0)--(6.5,2)--(8,2)--(8,0);
	\draw[-] (6.35,0) to (9.35,0);
	\draw[-] (6.5,-0.35)--(6.5,2.5);
    \draw[-, ultra thick] (6.5,2)--(9,0);
    \draw (8,1) node[right]{$_{\mathbb{E}(\chi)}$};
    \draw (7.25,2.1) node[above]{$_{\chi(\omega)}$};
    \draw (6.5,2) node{$_{-}$}; 
    \draw (6.5,2) node[left] {$_{1}$}; 
    \draw (6.5,-0.22) node[left] {$_{0}$};
    \draw (8,0) node{$_{|}$}; 
    \draw (8,-0.27) node[below] {$_{\omega}$};
    \draw (9,0) node{$_{|}$}; 
    \draw (9,-0.2) node[below] {$_{1}$};
    \draw [ultra thick,black] (8,0)--(9,0);
    \draw [ultra thick,black] (6.5,2)--(8,2);
\end{tikzpicture}
%\caption{Representación gráfica del valor esperado de la variable aleatoria $\chi:(0,1) \to L^{2}(0,1)$.}
\end{figure}

Although Problem \ref{Problem2} could have been solved directly by characterizing the set
$$
\left\{\omega\in(0,1)\,:\,\int_{(0,1)}\mathbf{1}_{(0,\omega)}(t)\,d\lambda(t)<\varepsilon^{p}\right\}=(0,\varepsilon^{p}),
$$
and applying the Lebesgue measure, the approach adopted here offers a much richer perspective. Inspired by the foundations of classical probability theory on $\mathbb{R}$, we have extended its essential notions —random variable, expectation, and integrability— to the setting of Lebesgue spaces $L^{p}$.

This generalization is not merely technical: it provides the starting point for the study of stochastic phenomena in Banach and Hilbert spaces, which naturally arise in advanced theories such as stochastic partial differential equations and Malliavin calculus; cf. \cite{Bogachev, DaPrato2014, DaPrato2004}.

The spaces $L^{p}(\Omega)$, with $p\in(1,\infty)$, are particularly well suited for this purpose. Their reflexivity (see \cite[Theorem 4.10]{Brezis}), a direct consequence of the Riesz representation theorem, ensures that every continuous linear functional can be identified with an element of the dual space $L^{q}(\Omega)$. This property allowed us to construct, in a natural way, the Bochner integral and the concept of an $L^{p}$-valued random variable, relying on the separable version of the Hahn–Banach theorem and on Pettis’s measurability theorem.

Taken together, these results reveal how the principles of measure theory and functional analysis converge in a single unifying idea: the understanding of probability as a form of linear representation within the Lebesgue spaces.

\section{Final Remarks}

The ideas developed throughout this work do not essentially depend on the specific properties of the Lebesgue spaces $L^{p}(\Omega)$, but rather on the fact that they constitute particular examples of separable Banach spaces. Indeed, both the elementary definitions and the arguments used in the proof of Pettis’s measurability theorem extend naturally to the case of functions taking values in an infinite-dimensional separable Banach space $V$ (see \cite[Chapter 2]{Diestel}).  

The separability of $V$ plays a crucial role: it allows the use of the countable version of the Hahn--Banach theorem, which ensures that the norm of any element in $V$ can be expressed as the supremum of a sequence of continuous linear functionals. This property makes it possible to describe the topology of $V$ through a countable family of continuous mappings and, therefore, to translate the notion of measurability into the setting of real-valued functions —the foundation of Pettis’s characterization.  

On the other hand, the completeness of the probability space $(\Omega,\mathcal{F},\mathbb{P})$ is equally fundamental. This condition guarantees that all subsets of null events belong to the $\sigma$--algebra $\mathcal{F}$, allowing one to work unambiguously with properties that hold almost surely. In particular, if the probability space were not complete, there could exist limit functions of sequences of simple random variables that are not measurable (see, for example, \cite{Durrett2019, JacodProtter2004}).  

It is worth mentioning that more general versions of Pettis’s measurability theorem exist which dispense with both the separability of the Banach space and the completeness of the probability space. However, in such cases one must assume that the image of the function $f:\Omega\to V$ is separable outside a null set, so that the essential behavior of the function can still be described within a separable subspace of $V$ (see \cite[Theorem 2, Chapter 2, Section 1]{Diestel}).  

Taken together, these observations highlight that the concepts introduced here —measurability, the Bochner integral, and expectation— fit within a broader theoretical structure in which probability intertwines naturally with functional analysis. The approach adopted, centered on Lebesgue spaces, sought precisely to provide a clear and concrete path into this interaction: through them, one can appreciate how the principles of measure theory and the linear structure of normed spaces merge to shape the fundamental notions of modern probability in abstract contexts.

\bibliographystyle{amsplain}

\end{document}